\documentclass{amsart}
\usepackage{amsmath,amssymb} 
\usepackage{amscd} 
\font\goth=eusm10

\newcommand\E{\mathcal E}

\newcommand\Ii{\hbox{\goth I}} 
 
\newcommand\Zz{\mathbf{Z}} 
\newcommand\ZZ{\mathbf{Z}}

\newcommand\Oc{\hbox{\goth O}}

\newcommand\Pq{\mathbf{P}^4} 
\newcommand\Pc{\mathbf{P}^5} 
\newcommand\Pj{\mathbf{P}} 
\newcommand\X{V_{2g-2}} 
 
\newtheorem*{RR}{Hirzebruch-Riemann-Roch theorem \normalfont{(\cite{Fulton} pg.290)}} 
\newtheorem*{CaMum}{Castelnuovo-Mumford Criterion 
\normalfont{(\cite{Mumford} pg.100)}} 
\newtheorem*{main}{Main Theorem}

\newtheorem*{proposition}{Proposition}

\newtheorem*{lemma}{Lemma} 
\begin{document} 
\title{ }
\centerline{\LARGE ACM vector bundles on prime Fano threefolds}
\medskip 
\centerline{\LARGE and complete intersection Calabi Yau threefolds} 
\bigskip  
\centerline{C.G.~Madonna}
\address{C.G.~Madonna \\
Dottorato di Ricerca in Matematica, Univerist\`a di Roma "Tor 
Vergata", Viale della Ricerca Scientifica, 00133 Roma, Italy} 
\email{madonna@mat.uniroma2.it, madonna@mat.uniroma3.it} 
\begin{abstract}
In this paper we derive a list of all the possible indecomposable 
normalized rank--two vector bundles without intermediate 
cohomology on the prime Fano threefolds and on the complete 
intersection Calabi-Yau threefolds, say $V$, of Picard number $\rho=1$.
For any such bundle $\E$, if it exists, we find the projective 
invariants of the curves $C \subset V$ which are the zero--locus
of general global sections 
of $\E$. In turn, a curve $C \subset V$ with such invariants is 
a section of a bundle $\E$ from our lists. This way we reduce the 
problem for existence of such bundles on $V$ to the problem for 
existence of curves with prescribed properties contained
in $V$. In part of the 
cases in our lists the existence of such curves on the general $V$ 
is known, and we state the question about the existence on the general 
$V$ of any type of curves from the lists. 
\end{abstract}
\maketitle

\medskip

\centerline{\bf \S \ 0. Introduction}

\medskip

Let $V$ be a smooth projective threefold, 
and let $\E$ be a vector bundle
on $V$. We say that $\E$ has {\it no intermediate
cohomology} if $h^i(V,\E \otimes L) = 0$,
$i = 1,2$ for any line bundle $L$ on $V$.
We say that an 
indecomposable (i.e. not split as a direct sum of two line bundles)
vector bundle on $V$ with no intermediate
cohomology is a {\it Cohen Macaulay bundle} (a CM bundle, for short). 
In the rank--two case we will call such a bundle {\it arithmetically
Cohen Macaulay} (an ACM bundle, for short).
In particular, let $\rho(V)= 1$ where $\rho(V)$ is the the rank of
the Picard group of $V$,
and let $D$ be the class of an ample generator of ${\rm Pic}(V)$
over $\ZZ$, i.e. ${\rm Pic}(V) \cong \ZZ[D]$.
For the line bundle $L = nD$ on $V$ denote by $\E(n)$
the twist of $\E$ by $nD$, i.e. $\E(n) = \E \otimes nD$.
Now, if the rank of $\E$ is two, 
the identity $\E^{\vee} = \E(-c_1)$, where $\det \E = c_1D$, 
together with the Serre duality, imply that $\E$ has no intermediate
cohomology if and only if $h^1(V,\E(n)) = 0$ for any $n \in \ZZ$. 
In this paper we will deal on ACM bundles on prime Fano threefolds
and complete intersection Calabi Yau threefolds. To start with,
by the criterion of Horrocks
the rank--two vector bundle $\E$ on $V = \Pj^3$
has no intermediate cohomology if and only it
splits as a direct sum of two line bundles
(see \cite{Hor} or \cite{OSS} -- Ch.2, Thm. 2.3.1), and hence
there are no ACM bundles on $\Pj^3$.
For the smooth quadric threefold $Q^3_2 \subset {\bf P}^4$ 
it was proved by Ottaviani that there exists a unique 
(up to twist) ACM bundle (see \cite{Ott}). 

In addition, by a result of Buchweitz, Greuel and Schreyer,
on any smooth hypersurface $X_d \subset \Pj^N$ of degree $d \ge 3$ 
(i.e. on any smooth $X_d$ different from $X_1 = \Pj^3$ or $X_2 = Q^3_2$) 
there exist infinitely many (up to twist) isomorphism classes
of CM bundles
(see \cite{BGS}). 

In \cite{Ma1} the author has obtained a cohomological criterion
for the splitting of a rank--two vector bundle $\E$ without intermediate 
cohomology on a smooth hypersurface $X_d$ in $\Pj^4$. 
It turns out that the approach used in \cite{Ma1} can be applied 
also to any smooth projective threefold $V$ for which the following condition
is satisfied:\par
\centerline{{\bf (${\ast}$)}
${\rm Pic}(V) = {\ZZ}[D]$ and $h^1(\Oc_V(nD))=0$ for any $n \in \ZZ$.} \par 
Moreover if the bundle $\E$ is {\it normalized} ,
i.e. $h^0(V,\E) \not= 0$ but $h^0(\E(-1)) := h^0(V,\E(-D))= 0$,
for $D$ being the ample generator of ${\rm Pic}(V)$,
then the zero--locus $C$ of a general global section of $\E$ is 
an 1-dimensional
(see \cite{H1} Remark 1.0.1) 
subcanonical curve on $V$ of well-defined genus $p$ and degree $d$, 
which will be dnoted by $C^p_d$.
In turn, the Serre construction (see e.g \cite{Bo})
makes it possible from a given curve $C \subset V$,
with the properties predicted after \cite{Ma1}
(provided such a curve on $V$ does exist),  
to construct a rank--two vector bundle $\E$ on $V$
without intermediate cohomology and such that $C$
is the zero--locus of a global section of $\E$.
This way, the problem for existence of
ACM bundles 
on the threefold $V$ with the property {\bf (${\ast}$)}
is reduced to the problem for existence of curves
$C \subset V$ with prescribed properties. 
This reduction has been successfully applied
by Arrondo and Costa (see \cite{AC}) to classify all ACM bundles
on the Fano threefolds
$Y_d \subset {\bf P}^{d+1}, d = 3,4,5$ of index $2$, 
as well in \cite{Ma2} to classify all such bundles on
the quartic threefold. We quote also the result by Szurek and Wi\'sniewski
in \cite{SW} in which some moduli spaces of rank--two bundles
on $Y_d$ are described by means of this reduction.

All the above gives us a ground to state the following
problem (depending on the choice of the threefold $V$): 

\medskip 

{\bf (${\ast}{\ast}$)} {\bf Problem}.
{\sl 
For any smooth projective threefold $V$ for which takes
place the condition {\bf (${\ast}$)},  
describe all the ACM bundles
on $V$. 
}

\medskip 

In this paper we study the problem {\bf (${\ast}{\ast}$)}
for two important
classes of projective threefolds -- the prime Fano threefolds
$V_{2g-2}, 2 \le g \le 12, g \not= 11$
and the Calabi Yau threefolds $X_r$, $r=5,8,9,12,16$, 
complete intersections in projective spaces.

\medskip 

By the well-known definition, 
the smooth projective variety $V$ of dimension $3$ 
is a Fano threefold if the anticanonical class $-K_V$ is ample. 
The Fano threefolds with Picard number $\rho=1$ 
are completely classified after the works of Iskovskikh
(see e.g. \cite{fanos}). 
Let $V$ be a Fano threefold with $\rho=1$, 
and let $D$ be the ample generator of ${\rm Pic}(V)$.
The positive integer $d = d(V) = D^3$ is called
the {\it degree} of $V$. 
Since ${\rm Pic}(V) = {\ZZ}[D]$ then 
$-K_V = rD$ for some positive integer $r = r(V)$
called the {\it index} of $V$.
Clearly any Fano threefold with $\rho= 1$
fullfills the condition {\bf (${\ast}$)}.  
By the classification of Fano threefolds with $\rho= 1$, 
the only Fano threefolds of index $r \ge 3$
are the projective 3-space $\Pj^3$ (for which $r = 4$)
and the 3-dimensional quadric $Q^3_2$ (for which $r = 3$).
The results of Horrocks for $\Pj^3$ in \cite{Hor}
and of Ottaviani for $Q^3_2$ in \cite{Ott} yield, of course,
the solution of the problem {\bf (${\ast}{\ast}$)} for Fano threefolds
of index $\ge 3$,
while the results of Arrondo and Costa in \cite{AC}
give the solution of {\bf (${\ast}{\ast}$)}
for the Fano threefolds $Y_d$ of index $2$ for which
the generator $H = (-K_Y)/2$ of ${\rm Pic}(Y)$
is very ample -- the cubic threefold $Y_3 \subset \Pj^4$,
the intersection of two quadrics $Y_4 \subset \Pj^5$
and the del Pezzo threefold $Y_5 \subset \Pj^6$. 

The Fano threefolds $V$ with $\rho= 1$ and
index $r = 1$ are called {\it prime}.
Obviously, the prime Fano threefolds are exactly
these Fano threefolds $V$ for which the anticanonical
class $-K_V$ is the ample generator of ${\rm Pic}(V)$
over $\ZZ$, i.e. ${\rm Pic}(V) = {\ZZ}[-K_V]$. 
The degree $d = (-K_V)^3$ of a prime Fano threefold $V$
is always even, and the integer $g = g(V) = \frac{d}{2}+1$
is called the {\it genus} of $V= V_{2g-2}$.
By \cite{fanos} prime Fano threefolds $V_{2g-2}$ exist
if and only if  
$2 \le g \le 12, g \not = 11$.

The main result of the paper is:

\begin{main} 
Let $\E$ be a normalized ACM bundle
on a prime Fano threefolds $V_{2g-2}$ of genus $g$.
Then $\E$ is a twist of one of the bundles in the list below:
\begin{enumerate} 
\item $c_1=-1$, $c_2=1$ and $\E$ is associated to a line 
in $V$;\par\noindent 
\item $c_1=0$, $c_2=2$ and $\E$ is associated to a conic 
in $V$;\par\noindent 
\item $c_1=1$, either $c_2=g+2$ and $\E$ is associated to a 
non-degenerate 
elliptic curve $C_{g+2}^1$ of degree $g+2$
contained in $V$ or $c_2=d <g+2$ and $C_{d}^{1}$ 
is degenerate;\par\noindent 
\item $c_1=2$, $c_2=2+2g$ and $\E$ is associated to a 
curve $C_{2g+2}^{g+2}$ of genus $g+2$ and degree $2g+2$ contained in $V$; 
\par\noindent 
\item $c_1=3$, $c_2=5g-1$ and $\E$ is associated to a smooth 
$2$-canonical curve $C_{5g-1}^{5g}$ contained in $V$. 
\end{enumerate} 
\end{main}

In showing previous theorem,
our starting point is that 
the procedure of \cite{Ma2} and \cite{Ma1}
works word by word for any prime Fano threefold $V_{2g-2}$
for which the anticanonical class 
determines an embedding in $\Pj^{g+1}$.

The only prime Fano threefolds, 
for which the anticanonical class does not define an embedding
in a projective space, are the sextic double solid 
and the double quadric (see \cite{fanos}).  
To include also these last in our main result
we first give a generalization of the splitting
criteria in \cite{Ma1}.

\medskip

We remark also 
the same approach works for
the case of complete intersection Calabi Yau threefolds.
By definition a Calabi Yau (CY, for short) threefold is a smooth 
3-dimensional projective variety $X$ such that the canonical 
class $K_X = 0$. In particular, if 
$X = X_{r_1,...,r_k} \subset \Pj^{k+3}$ 
is a complete intersection of hypersurfaces  
of degrees $r_1,...,r_k$ then $X$ is called a complete 
intersection CY threefold (ciCY, for short). Clearly the product 
$r = r_1 \cdot ... \cdot r_k$ is the degree of $X \subset \Pj^{k+3}$. 
There exist five types of ciCY:
the quintic $X_5 \subset \Pj^4$, 
the intersection $X_9 \subset \Pj^5$ of two cubics,
the intersection of a quartic and a quadric $X_8$ in $\Pj^5$,
the intersection of a cubic and two quadrics $X_{12}$ in $\Pj^6$,
and the intersection of 4 quadrics $X_{16}$ in $\Pj^7$.\par
Also for these threefolds we will give in the last section
a list of eventually all normalized ACM bundles. We will omit 
its proof, being a very straightforward computation of 
the Fano's case.

\medskip

Both for Fano's and ciCY's we
reduce the problem of existence
of ACM bundles 
to the problem of existence of curves with prescribed
properties on them.
We quote, of course lines and conics exist on all the threefolds
from the two lists, but it is a problem to
prove the existence of all the types of curves predicted 
by these two lists.

We have also to point out that in the recent paper
\cite{Bo} of Beauville the pfaffian representation of
a hypersurface $Z_r$ in ${\bf P}^N$, $N=3,4$, of degree $r$ is related
to an ACM bundle $\E$ on $Z_r$, and the existence of
such a bundle (i.e. of a pfaffian representation of $Z_r$) is
showed to be 
equivalent to the existence of special codimenion 2 cycles $C^{pf}$ 
on $Z_r$. Although the numerical estimates 
imply immediately a list of all the expected $(r,N)$ for which
the general $Z_r$ can't have a pfaffian representation,
the problem of existence such cycles $C^{pf}$
(hence -- of pfaffian representations) on the general
$Z_r$ in the remaining cases had been solved completely
by involving computer programs (see the appendix of
F.O. Schreyer in \cite{Bo}).
In certain cases the cycles $C^{pf}$ have been
constructed directly on the base of the geometry of the variety
-- see e.g. \cite{MT} and \cite{IM} for $(r,N) = (3,4)$
and \cite{IM} for $(r,N) = (4,4)$. 

Because the problem of existence of all the types of curves 
predicted
on both the prime Fano threefolds and the 
ciCY threefolds 
remains open, we state however the following:

\medskip

\noindent {\bf Conjecture.}
{\sl
All the types of curves (hence -- of rank--two bundles)
from the two lists 
exist.
}

\medskip 

In connection to the results in \cite{Bo},
in \S 2.2 below we show that the vector bundles
from the list with $c_1 = 2$ and $c_2 = 10$
on the Fano threefold $V_6$ are parameterized by the
pfaffian cubic 4-folds containing $V_6$. 
In addition, we show the existence and
describe some of bundles from the list 
for prime Fano threefolds of small genera -- see \S 2. 

\medskip

{\bf Acknowledgment.} I wish to thank Professor A.\ Iliev for many
useful conversations during his visit in Mathematics Department
of Roma Tre on Genuary 2000. 

\medskip

\centerline{\bf \S \ 1. Notation and generalities} 

\medskip 

In this paper we work  over the complex field $\mathbb{C}$.
Our main references are \cite{H},\cite{fanos}, and \cite{OSS}.
An $n$-space is a projective
space of dimension $n$. 
A curve contained in a projective threefold $V$ 
is a reduced (eventually reducible) 
1-cycle contained in $V$, hence an element of the Chow
group $A^2(V)$. 

When $\text{Pic}(V) \cong \mathbb{Z}$  
we use this isomorphism to identify line bundles with 
integers; in particular, for any  vector bundle $\E$ we denote by
$c_1(\E)=c_1 \in \mathbb{Z}$. 
A locally complete intersection curve $Y \subset X$ is called 
$a$-subcanonical if the canonical class
$\omega_Y=\Oc_Y(a)$ for some $a \in 
\mathbb{Z}$.\par 
We introduce
$b(\E)=b=\max \{ n \mid h^0(\E(-n)) \ne 0 \}$, hence
$\E$ is
normalized
if and only if $b(\E)=0$.
Note that 
any rank--two vector bundle may be normalized by mean of the twist $\E(-b)$,
so we may always assume the bundle is normalized.
We use the notion of stability given in \cite{OSS}; 
when $\text{Pic}(V) \cong \mathbb{Z}$ (as for prime
Fano and ciCY threefolds), using our notation, a rank--two 
vector bundle $\E$ is stable if and only if $2b-c_1<0$. 
Note also that the number $2b-c_1$ is invariant by twisting i.e 
$2b(\E)-c_1(\E)=2b(\E(n))-c_1(\E(n)) \forall n \in \mathbb{Z}$.
Let us denote by $V$ a prime Fano threefold or a ciCY threefold, and
let $H$ denote the class of the ample generator of 
${\rm Pic}(V)$ 
and by $l$ the class of a line on $V$. 
Then $H^3=\deg V$, $H \cdot l=1$ 
and $H^2=(\deg V)l$.
If $\E$ is a rank--two vector bundle on $V$ then for any $n \in \Zz$ 
we have 
\[ 
\begin{array}{ll} 
c_1(\E(nH)=c_1(\E)+2nH=(c_1+2n)H \\ 
c_2(\E(nH))=c_2(\E)+c_1(\E) \cdot nH +(nH)^2= 
\left( c_2+(\deg V)nc_1+(\deg V)n^2 \right)l. \\ 
\end{array} 
\] 
We will frequently use the following:

\begin{RR} 
The following holds:\par\noindent 
(i) $\chi(\Oc_{\X}(m))=\frac16(g-1)(m+1)(2m+1)m+2m+1$ for any $m \in 
\Zz$;\par\noindent 
(ii) $\chi(\E)=(2g-2)c_1^2 \left( \frac{c_1}{6}+\frac14 \right)-\frac{c_2}2(c_1+1)+\frac{c_1}6 
(g+11)+2$, where  
$\E$ is a rank--two vector bundle on $\X$ with Chern classes $c_1(\E)=c_1$ and $c_2(\E)=c_2$. \par\noindent
(iii)
$\chi(\E)=\frac{r}{6}c_1^3-\frac{c_1c_2}{2}+\frac{c_1}{12}[12(k+4)-2r)]$
where $\E$ is a rank--two vector bundle on $X_r \subset
\Pj^{k+3}$ with Chern
classes
$c_1$ and $c_2$;
\par\noindent
(iv) $\chi(\Oc_{X_r}(nH))=\frac{r}{6}(n^3+5n)$.
\end{RR} 

\begin{CaMum} 
Let $V$ be a non-singular projective variety 
with ${\rm Pic}(V) \cong \mathbf{Z}$, and 
let $\E$ be an $m$-regular coherent sheaf on $V$ 
i.e. such that $h^i(\E(m-i))=0$ for all $i>0$. 
Then $\E(k)$ is globally generated 
if $k \geq m$. 
\end{CaMum} 

\medskip

\centerline{\bf \S \ 2. A proof of the Main Theorem.}

\medskip

\begin{lemma}
Let $V_{2g-2}$ be a prime Fano threefold and let $D$ be the ample generator
of ${\rm Pic}(V_{2g-2})$. Then for any curve $C \subset V_{2g-2}$ there
exists
a 1-cycle $L=D_1 \cdot D_2$, $D_1,D_2 \in \mid D \mid$
which does not intersect $C$.
\end{lemma}

We can now give the following generalization of the
splitting criteria in \cite{Ma1}. By previous lemma, its
proof, is a very
straightforward procedure of the one as in \cite{Ma1}.

\begin{proposition}
Let $V_{2g-2}$ be a prime Fano threefold and $D$ the class of
an ample generator of $V_{2g-2}$.
Let $\E$ be a normalized
rank--two vector bundle on $V_{2g-2}$. Then: if $\oplus_n H^1(nD)=0$,
$\E$ splits as a direct sum of two line bundles, unless \par
\centerline{
$-2L \cdot D -L \cdot K_L \leq -L \cdot c_1(\E) \leq L \cdot K_L$.}
\end{proposition}

\par\noindent
{\bf Proof of the Main Theorem.} Let us start as in 
\cite{Ma2} for the case $g=3$, with a case 
by case analysis.
By previous proposition we find $-4<-c_1<2$.\par\noindent
If $c_1=-1$ we find that $\chi(\E(-1))=-g-2=-3-g+c_2$ hence $c_2=1$ and 
$\E$ is associated 
to a line on $V$. \par \noindent 
If $c_1=0$ we find $\chi(\E)=-\frac{c_2}{2}+2=1$ hence $c_2=2$ and 
$\E$ is associated to a 
conic.\par \noindent 
If $c_1=1$ since $2p-2=0$, $\E$ corresponds to an elliptic curve 
of degree $c_2$. Since 
\[ 
h^0(\E)=\chi(\E)=(2g-2)\left( \frac16+\frac14 \right) 
-c_2+\frac{g+11}{6}+2=1+h^0(\Ii_C(1)) 
\] 
it follows that $c_2=g+2-h^0(\Ii_C(1))$. 
If $C$ is not contained in any hyperplane then $h^0(\Ii_C(1))=0$ and 
the degree of $C$ is $g+2$. 
In any case $c_2 \geq 3$. 
\par \noindent 
If $c_1=2$ since $\chi(\E(-1))=0$ we find that $\E$ is associated to 
a curve of 
genus $p=g+2$ and degree $d=2+2g$. 

\par

\noindent
If $c_1=3$ since $h^3(\E(-3))=0$ by Castelnuovo-Mumford regularity
a general global section of $\E$
is a smooth curve of genus $p$ and degree $p-1$.
Since $\chi(\E(-1))=0$ then $c_2=5g-1$ and $p=5g$. 

\medskip

\centerline{\bf \S \ 3. Examples.}

\medskip

{\bf \S \ 3.1. Elliptic cubics on $V_6$.} Any prime 
Fano threefold $V=V_6 = Q \cap F \subset \Pj^5$
of genus $g=4$ is a complete intersection
of a quadric $Q$ and a cubic hypersurface $F$ in $\Pj^5$
(see \cite{fanos}).
The quadric $Q$ is uniquely defined by $V$,
and if $V$ is smooth then $Q$ has rank
5 or 6.
For, if ${\rm rk}(Q) = r \le 4$ then the vertex
$\Pj^{6-r-1} = {\rm Sing}(Q)$ of $Q$ will be
a space of dimension $6-r-1$, and then
$\Pj^{6-r-1} \cap F \subset {\rm Sing}( V) $ will
be non-empty.
If ${\rm rk}(Q) = 6$ then $Q \cong Gr(2,4)$
is isomorphic to the Grassmannian of lines
in $\Pj^3$ embedded by the Pl\"ucker map in $\Pj^5$.
Then on $Q$ lie two $\Pj^3$--families
$\Lambda = \{ \Pj^2(t); t \in \Pj^3 \}$
and
$\overline{\Lambda} = \{ \overline{\Pj}^2(t); t \in \Pj^3 \}$
of planes. 

For a general $F$, the prime Fano threefold $V$
does not contain planes, 
and hence the cubic $F$ intersects on any plane $\Pj^2(t) \subset Q$
(or $\overline{\Pj}^2(t) \subset Q$) 
a 1-cycle of degree $3$ and it is easy to see that
for the {\it general} $V$ the 1-cycles is
a smooth plane cubic on $V$.
Let, in turn, $C \subset V$ be a plane cubic,
i.e. $C = V \cap \Pj^2$ for some plane $\Pj^2 \subset Q$,
and assume for simplicity that $C$ is smooth.
The argument below, which depends only on the degree
of the 1-cycle $C = V \cap \Pj^2$ takes place also
for reducible and non-reduced intersections $V \cap \Pj^2$.

Since ${\rm deg}(C) = 3$ and $C \subset V \subset Q$
then the plane $\Pj^2 = Span(C)$ will lie on $Q$;
therefore the plane $\Pj^2$ will belong to one of the
two families $\Lambda$ or $\overline{\Lambda}$.

If ${\rm rk}(Q) = 5$ then $Q$ is a cone
with vertex a point,
over a smooth $3$-dimensional quadric.
Then on $Q$ lies one $\Pj^3$--family
$\Lambda = \{ \Pj^2(t); t \in \Pj^3 \}$ of 2--spaces.  
The general plane $\Pj^2$ will intersect on $V$
a smooth plane cubic curve. 
The same argument as above implies that
the plane cubics on $V$ are, in fact, the intersections
$C = V \cap \Pj^2(t)$ for $\Pj^2(t) \in \Lambda$.
This yields the following

\begin{proposition}
Let $V_6 = Q \cap F \subset \Pj^5$ be a prime Fano threefold 
of genus $4$. Then:

(1) if ${\rm rk}(Q) = 6$ then on $V_6$ lie two $\Pj^3$--families
of plane cubic intersected by the elements of the
two $\Pj^3$--families of planes on $Q$; moreover if $V_6$ is
general then the general element of any of these two families
is a smooth plane cubic curve;

(2) if ${\rm rk}(Q) = 5$ then on $V_6$ lies one $\Pj^3$--family
of plane cubics intersected by the elements of the
unique $\Pj^3$--family of planes on $Q$; moreover if $V_6$ is
general then the general plane cubic on $V_6$ 
is a smooth plane cubic curve.
\end{proposition}

\medskip

{\bf \S \ 3.2. Elliptic quartics on $V_6$.} On the smooth 
Fano threefold $V=V_6 = Q \cap F \subset \Pj^5$ lie
a 2-dimensional family of conics
$\Gamma = \{ q : q \subset V, \ deg(q) = 2 \}$
(i.e. any component of $\Gamma$ has dimension 2);  
and the general conic $q \subset V$
(i.e. the general element of any component of $\Gamma$) 
is smooth (see \cite{fanos} for example).

Let $q$ be a smooth conic on $V$, and let
$\Pj^3 \subset \Pj^5$ be a general 3-space which
contains $q$.
Since the smooth threefold $V$ does not contain 2-dimensional
quadrics then the 1-cycle
$C^4_6 = V \cap \Pj^3 = Q \cap F \cap \Pj^3 = q + C$
will be a complete intersection of a quadric and a cubic
in $\Pj^3$ containing $q$ as a component.
Assuming that $V$ is {\it general}
it is easy to see that for the general
$\Pj^3 \supset q$ the residue 1-cycle $C = C^1_4$
is a smooth elliptic quartic on $V$.  
Let, in turn, $C = C^1_4 \subset V$ be an
elliptic quartic. Clearly $C$ can't lie on a plane,
otherwise, for reason of intersection number,
this plane will lie also on $V$.
Therefore $C$ is a complete intersection
of two quadrics in $\Pj^3 = Span(C) \subset \Pj^5$.
Since the quadric $Q \subset \Pj^5$
has rank $\ge 5$ then $Q$ does not contain 3-spaces,
and $S_2 = Q \cap \Pj^3$ will be a quadric surface
containing $C$. Similarly $S_3 = F \cap \Pj^3$
is a cubic surface which contains $C$, and which
does not have components which lie on $V$.
Therefore the 1-cycle 
$C^4_6 = C + q = S_2 \cap S_3 = Q \cap F \cap \Pj^3 \subset Q \cap F = V$
is a cubic section of $S_2$, containing $C = C^1_4$ as a component.
Since $C$ is a quadratic section of $S_2$, then the residue
$1$-cycle $q$ is a linear section of $S_2$, i.e.
$q$ is a conic on $V$.
This yields the following

\begin{proposition}
Let $V_6 = Q \cap F \subset \Pj^5$ be a prime Fano threefold 
of genus $4$. Then on $V_6$ there exist elliptic quartics
(i.e. effective 1-cycles which are complete intersections
of two quadrics in a $\Pj^3$); and if $ V_6$ is general
then the general elliptic quartic
on $V_6$ is a smooth. 
\end{proposition}

\medskip

{\bf \S \ 3.2. Elliptic quartics on $V_8$.} Any 
prime Fano threefold $V_8$ of genus $g=5$,
is a complete intersection of three quadrics in ${\bf P}^6$
(see \cite{fanos}). 
Let us denote by ${\Pj}^2[N]$
the plane of quadrics 
whose base locus defines $V_8$. 
The general quadric $Q \in {\Pj}^2$ is smooth (i.e. of rank = $7$), 
and it is well known, if the quadrics of the net are not all
singular, the
Hessian 
$N_o = \{ Q \in {\bf P}^2[N]: {\rm rk}(Q) \le 6 \}$
is a plane curve of degree 7.

For the general $V_8$ the curve $N_o$ is smooth,
and any any element $Q \in N_o$ is a cone, with vertex a point, 
over a smooth
5--dimensional quadric.  
Therefore any quadric $Q \in N_o$ contains two ${\bf P}^3$--families
${\Lambda}_Q$ and $\overline{\Lambda}_Q$ of projective $3$--spaces,
and since any of these two families defines uniquely
the quadric $Q$, there exists a $2$--sheeted unbranched
covering
${\pi} : \tilde{N}_o \rightarrow N_o$,
${\pi}^{-1}(Q) = \{ {\Lambda}_Q , \overline{\Lambda}_Q \}$,
defined by a $2$-torsion sheaf $\nu \in {\rm Pic}_2(N_o)$
(see \cite{Bo-Prym}).  
First we prove the following

\begin{lemma}
There are no smooth plane curves of degree $d \geq 3$ 
on the prime Fano threefold $V_8 \subset \Pj^6$.  
\end{lemma} 
 
\begin{proof} 
Let $C \subset V_8=V$ be a smooth plane curve of
degree $d \geq 3$, and let $\Pj^2$ the plane containing it.
Clearly, for reason of intersection number, the plane
$\Pj^2$ can't lie in $V$.
Since $\deg(C) \ge 3$ then the general line 
$L \subset \Pj^2$ intersects $C \subset V$
in at least 3 points. But then $L$ intersects also
$V \supset C$ in these points and hence
it is contained in any quadric containing $V$, and hence in $V$.
Since $L \subset \Pj^2$ is general
then $L \subset V$ implies that $\Pj^2 \subset V$, which is absurd.
\end{proof} 
 
\begin{proposition} \label{prop:8,(1,4)} 
On the general prime Fano threefold 
$V_8 \subset \Pj^6$ 
there exists a projectively normal elliptic 
curve of degree $4$; moreover, any elliptic
quartic on $V_8$ is intersected by a $3$-space which lies
in a rank $6$ quadric containing $V_8$. 
\end{proposition} 
 
\begin{proof}
Let $Q \in N_o$ be a rank $6$ quadric containing $V_8=V$.
Then on $Q$ lie two $\Pj^3$ families of $3$-spaces;
and let $\Pj^3$ be one of these $3$-spaces.
Since $V$ is a complete intersection of
$Q$ and two quadrics, say $V = Q \cap Q' \cap Q''$,
and since $\Pj^3 \subset Q$,
then $V \cap \Pj^3 = Q' \cap Q' \cap \Pj^3$.

Since $Q'$ and $Q''$ are quadrics
then any of the two cycles
$S' = Q' \cap \Pj^3$ and $S'' = Q'' \cap \Pj^3$
is either $\Pj^3$ or a 2-dimensional quadric.

But if some of these two cycles, say $S''$, is $\Pj^3$,
then $V$ will contain the cycle
$Q' \cap (Q'' \cap \Pj^3) = Q' \cap \Pj^3 = S'$.
Therefore
either $V \supset S''= \Pj^3$, and then $V$ will be not irreducible,
or $S'' \subset V$ will be a quadric surface, and both cases
are impossible.

Therefore $S'$ and $S''$ are quadric surfaces in $\Pj^3$,
and $V$ contains the intersection cycle 
$C= Q' \cap Q' \cap \Pj^3$ =
$(Q' \cap \Pj^3) \cap (Q'' \cap \Pj^3) = S' \cap S'' \subset \Pj^3$.
As above $C$ has not 2--dimensional components, and hence
it is a 1--cycle.
It is easy to see that for the general $V$ the general
such $C$ is smooth, i.e. $C$ is a smooth elliptic quartic curve. 
Let, in turn, $C$ be an elliptic quartic on $V$,
i.e. $C = C^1_4 \subset V \subset \Pj^6$
is a complete intersection of two quadrics in the
3-space $\Pj^3 = Span(C)$.

Since $C \subset V$ and $V$ is a complete intersection
of three quadrics in $\Pj^6$, and since 
$C \subset V \cap \Pj^3$ is a complete intersection
of two quadrics in $\Pj^3 = Span(C)$,
then there exists a quadric $Q \subset \Pj^6$
containing $V$ such that $\Pj^3 \subset Q$. Since the rank 7 quadrics
in $\Pj^6$ do not contain 3--spaces, then $rank \ Q \le 6$,
i.e. $Q \in N_o$.
\end{proof} 

\medskip

{\bf \S \ 3.3. Canonical curves $C^6_{10}$ on $V_6$.} \par
Let us introduce the following:\par\noindent
{\bf Definition.}
{\it
Let ${\rm Pfaff}$ be the set of pfaffian cubic 4-folds 
in $\Pj^5$. Let us denote by $M=\Pj^6_F \cap {\rm Pfaff}$ 
the 5-dimensional subvariety of pfaffian  
cubic 4-folds containing $V_6$. 
We will call the 5--fold of pfaffian cubics
$\tilde{F}$ not of type $L \cdot Q$
for some linear form $L$, say 
$M_o$, the set of pfaffian representations of $V_6$. }
 
Next we will show the role 
played by pfaffian representations of the 
threefold $V_6$. Following \cite{Bo} 
we will enlight the analogy with the case 
of the quartic threefold ($g=3$).
 
\begin{proposition} Let $V_6$ be general. 
Then, there exist a bijective map between the family of 
pfaffian representations 
$M_o$ of $V_6$ and the family of curves of genus $6$
and degree $10$ on $V_6$.
\end{proposition} 
 
\begin{proof} 
Let us denote by $V = V_6=Q \cap F \subset \Pc$. 
We will show that the general element of $M_o$, 
i.e. the general pfaffian representation of $V$, corresponds 
to the general rank--two bundle $\E$ on $V$ associated 
by the Hartshorne--Serre correspondence 
to a curve $C = C^6_{10} \subset V$ 
of genus $6$ and degree $10$ contained in $V$. 
To show this we may proceed 
in the following way (or as in \cite{IM}). 
To an element in $M_o$ corresponds an unique pfaffian cubic 
4-fold $F \subset \Pc$ such that $F \cap Q=V$. Such a 4-fold 
contains 
(see  \cite{Bo} Proposition 8.2 
for a proof) 
a quintic Del Pezzo surface $S_5 \subset V$
and $S_5 \cap Q=C_{10}$ is a curve of degree $10$ 
contained in $V$. 
Since $K_{S_5}=\Oc_{S_5}(-1)$, by adjunction, 
$C_{10}=C$ is an arithmetically normal curve
of genus $6$ and degree $10$ contained in $V$.
The bundle associated to $C$ fits in the following exact sequence
$$
0 \to \Oc_V \to \E \to \Ii_C(2) \to 0,
$$
it does not split being $c_2=10$ and has no intermediate cohomology
being $C$ arithmetically normal. For this bundle we find, $b=0$,
and $c_1=2$.
\end{proof}

{\bf \S \ 3.4. The Gushel's approach to the classification
of prime Fano threefolds $V_{10}$ and $V_{14}$
(\cite{Gushel1} and \cite{Gushel2}).}
To study $V_{10}$
Gushel first proves that on any smooth prime Fano threefold
$V_{10}$ there exists a smooth curve of degree $4$ and genus $1$.
This curve is contained in a smooth hyperplane section $H$ of $V_{10}$
and he associates to the pair $(Z,H)$ a rank--two vector bundle $\E$
generated by global sections with $c_1=H$, $c_2=Z$, and $h^0(\E)=5$.

In our setting this means that there exists a rank--two
vector bundle $\E$ on $V_{10}$ with $c_1=1$, $c_2=4$, with $h^0(\E)=5$.
Since $\E$ is globally generated, a global section $s \in H^0(\E)$ 
vanishes along a smooth elliptic quartic $Z$ 
which is a complete intersection of two quadrics in $\Pj^3$.
Since $Z$ is arithmetically normal $\E$ has no intermediate cohomology,
and $\E$ does not split since $Z$ is not a complete intersection in $V_{10}$.
The bundle $\E$ determines a morphism $\psi:V_{10} \to Gr(2,5)$.
In \cite{Gushel1} Gushel proves that $\deg \ \psi$ = 1 or 2,
which implies the two types of Fano threefolds of genus $6$. 
Using the same procedure Gushel shows the existence of a smooth
elliptic quintic $Z$ on $V_{14}$ which implies the existence of a
globally generated rank--two vector bundle on $V_{14}$ with
Chern classes $c_1=1$ and $c_2=5$, and $h^0(\E)=6$.
Since the curve is projectively normal $\E$ has no intermediate cohomology;
also in this case $\E$ does not split since $Z$ is not a complete intersection
in $V_{14}$ for reason of the degrees.
The bundle $\E$ defines a morphism $\psi:V_{14} \to Gr(2,6)$, and
in \cite{Gushel2} Gushel proves that $\psi$ is an embedding. 

\medskip

\centerline{\bf \S \ 4. The case of ciCY's.}

\medskip

For ciCY threefolds we state the following:

\begin{proposition}
Let $X_r$ be a ciCY threefold,
and let $\E$ be a normalized ACM bundle on it.
Then $\E$ is a twist of one of the bundles in the list below:
\begin{enumerate}
\item $c_1=-2$, $c_2=1$ and $\E$ is associated to a line;
\item $c_1=-1$, $c_2=2$ and $\E$ is associated to a conic;
\item $c_1=0$, $3 \leq c_2 \leq \deg(X)$ and $\E$ is associated to an
elliptic curve $C^1_{c_2}$;
\item $c_1=1$, $c_2=2r-2$ and $\E$ is associated to a curve $C_{2r-2}^r$;
\item $c_1=2$, $c_2 \leq 3r-1$ and $\E$ is associated to a curve
of degree $c_2$ and genus $c_2+1 \leq 3r$; 
\item $c_1=3$, $c_2=4r$ and $\E$ is associated to a curve $C_{4r}^{6r+1}$;
\item $c_1=4$, $c_2=6r$ and $\E$ is associated 
to a smooth curve $C^{12r+1}_{6r}$.
\end{enumerate}
\end{proposition}

We omit the proof of previous proposition, being a very straightforward
computation of the technique used for the Fano's case.\par\noindent
{\bf Example.} If $X=X_5$ is a general quintic hypersurface
the existence of bundles in 1.--4. from the list follows by
the existence of the corresponding curves on the threefold.
The existence of bundles as in 7. of the above list was showed
in \cite{Bo}.

\end{document}